\documentclass{elsart}
\usepackage{amsmath,amssymb}

\newtheorem{theorem}{Theorem}
\newtheorem{corollary}{Corollary}
\newtheorem{definition}{Definition}

\begin{document}

\begin{frontmatter}

\title{Distribution of the full rank in residue classes for odd moduli}

\author{William J. Keith}

\address{Drexel University, 3141 Chestnut, Philadelphia, PA 19104, USA}
\ead{wjk26@drexel.edu}

\begin{abstract}
The distribution of values of the full ranks of marked Durfee symbols is examined in prime and nonprime arithmetic progressions.  The relative populations of different residues for the same modulus are determined: the primary result is that $k$-marked Durfee symbols of $n$ equally populate the residue classes $a$ and $b$ mod $2k+1$ if $gcd(a,2k+1)=gcd(b,2k+1)$.  These are used to construct a few congruences.  The general procedure is illustrated with a particular theorem on 4-marked symbols for multiples of 3.
\end{abstract}


\end{frontmatter}

\section{Results}

In 2007 Andrews \cite{Andrews} examined the moments of the rank function on partitions, describing a new combinatorial object called $k$-marked Durfee symbols and a statistic called the full rank which described interesting aspects of its behavior.  The full rank for 2-marked and 3-marked Durfee symbols is equidistributed in certain arithmetic progressions mod 5 and 7, giving rise to congruence theorems, and nearly equidistributed in all other progressions over those moduli.  The present author's doctoral thesis was partially devoted to explaining that behavior.  The examples given were shown to be the simplest instances of an infinite family of relations on all full ranks; while most of the other relations do not produce clean congruences their failure mode is fairly neat and is fully explicated.  This is a reprint and slight updating of that material for wider circulation.  The primary results are:

\begin{theorem}\label{Main Theorem} Let $c=2l+1 \in \mathbb{N}$.  Say $NF_l (j,c,cn+d)$ is the number of $l$-marked Durfee symbols of $cn+d$ with full rank congruent to $j$ \emph{mod} $c$.  Then, if $gcd(i,c) = gcd(j,c)$, we have $NF_l (i,c,cn+d) = NF_l (j,c,cn+d)$.
\end{theorem}

As a corollary,

\begin{corollary}\label{PrimeCorollary}
If $c=2l+1$ is prime, then $NF_l (i,c,cn+d) = NF_l (j,c,cn+d)$ for all $i,j \not\equiv 0 \ \emph{mod} \ c$.
\end{corollary}

This is the case for Andrews' two theorems previously discussed.  For those theorems complete equidistribution in residue classes comes about due to a second consequence that will be easily seen from Theorem \ref{Main Theorem}'s method of proof:

\begin{theorem}\label{PrimeLMinusOne}
If $c=2l+1$ is prime, then 
\[
NF_l (0,c,cn+d) - NF_l (1,c,cn+d) = N(l-1,c,cn+d) - N(l,c,cn+d) \text{.}
\]
\end{theorem}

\noindent where $N(b,c,n)$ is the number of standard partitions of $n$ with usual (Dyson's) rank congruent to $b$ modulo $c$.  Since the original purpose of the full rank was to investigate the behavior of the overall number of $l$-marked Durfee symbols of $n$, denoted $D_l (n)$, we can combine Corollary \ref{PrimeCorollary} and Theorem \ref{PrimeLMinusOne} to obtain  

\begin{corollary}\label{PrimeOverall}
If $c=2l+1$ is prime, then 
\[
\sum_{m=-\infty}^{\infty} NF_l (m,n) = D_l (n) \equiv N(l-1,c,n) - N(l,c,n) \, (mod \, c) \text{.}
\]
\end{corollary}

Because this difference is 0 for $c=5, d=1,4$ and $c=7, d=0,1,5$, we have full equidistribution and a clean congruence theorem in those progressions.

Some notes on concurrent work are in order.  Independently, Kathrin Bringmann, Frank Garvan, and Karl Mahlburg also studied the automorphic aspects of rank moments and Durfee symbols; between the original writing of this thesis and the current reprinting, they have published an article in INRM (\cite{BGM}).  While there is some overlap, their main interest is rank moments; this paper focuses more on Durfee symbols and the full rank.  Here we are more enumerative, establishing differences of the full rank in residue classes and determining congruences for those values in terms of the usual rank.  An example of a theorem of interest to us is 

\begin{theorem}\label{FourMarkedModThree} $D_4 (3n) = \eta_6 (3n) \equiv 0 \, (mod \, 3)$.
\end{theorem}

Bringmann and Ben Kane, in addition, are preparing a paper studying 2-marked Durfee symbols in general arithmetic progressions; that paper overlaps with this one for the modulus 5 case, for which we have both proved the same result independently.  Kathy Ji has a paper in the arXiv \cite{Ji} on several bijections for Durfee symbols and odd Durfee symbols, with their combinatorial implications.  Other doctoral students of George Andrews are also studying Durfee symbols, particularly Kagan Kursungoz (\cite{Kursun}).  That work is joint with Cilanne Boulet, who also has a paper in preparation studying their symmetries: \cite{Boulet}.  Boulet's previous work, on Garvan's $k$-rank (first as a doctoral student under Richard Stanley \cite{BouletThesis} and later with Igor Pak (\cite{BouletPak}) and individually (arXiv:math/0607138)) may also be of interest to the reader seeking information on the behavior of the standard rank and its generalizations.

\section{Definitions}

A vector $\lambda = ( \lambda_1 , \dots \lambda_k ) \in
\mathbb{N}^k$ is a \emph{partition} of $n$ if $\lambda_1 \geq
\dots \geq \lambda_k \geq 1$ and $\sum \lambda_i =
n$. The number of partitions of $n$ is denoted $p(n)$.  Dyson's
\emph{rank} of $\lambda$ is $\lambda_1 - k$.  Let $N(m,n)$
be the number of partitions of $n$ with rank $m$:

\begin{equation}\label{RankGenFn}
R_1(z;q) := \sum_{m=-\infty}^{\infty} \sum_{n \geq 0} N(m,n) z^m
q^n = \sum_{n \geq 0} \frac{q^{n^2}}{{(zq;q)}_n {(q/z;q)}_n} \, \text{,}
\end{equation}

\noindent where ${(a;q)}_n = \prod_{i=0}^{n-1} (1-aq^i)$.  It is
easily observed that $N(m,n) = N(-m,n)$.

Dyson was motivated by the fact that partitions with rank $\equiv
i \, (mod \, 5,7)$ are distributed evenly for partitions of $5n+4$
and $7n+5$ respectively: that is, if $N(i,c,n)$ denotes the number
of partitions of $n$ with rank $\equiv i \, (mod \, c)$, then
$N(i,5,5n+4)=N(j,5,5n+4)$ and $N(i,7,7n+5) = N(j,7,7n+5)$ for all
$i,j$.  This provided a combinatorial explanation of Ramanujan's
famous theorems that $p(5n+4) \equiv 0 \, (mod \, 5)$ and $p(7n+5)
\equiv 0 \, (mod \, 7)$.

In studying further partition congruences, A.O.L. Atkin and Frank
Garvan \cite{AtkinGarvan} constructed the $k$-th \emph{moments} of
the rank function.  George Andrews \cite{Andrews} has in turn
constructed the \emph{symmetrized k-th moments}

\[ \eta_k (n) = \sum_{m=-\infty}^{\infty} \left( {{m + \lfloor \frac{k-1}{2}} \rfloor \atop k} \right) N(m,n) \, \text{,} \]

\noindent and associated to these the \emph{k-marked Durfee
symbol}, in which two rows of parts are marked with $k$ subscripts or
colors, according to the following rules:

\begin{definition} The ordered, subscripted vector pair ${ \left( \begin{matrix} t_1 & t_2 & \dots & t_r \\
b_1 & b_2 & \dots & b_s
\end{matrix} \right) }_c $ is a $k$-marked Durfee symbol of $n=c^2 + t_1 + \dots t_r + b_1 + \dots + b_s$ if
\begin{itemize}
\item $t_i, b_j \in \{ 1_1, 1_2, \dots , 1_k, 2_1, 2_2, \dots ,
2_k, \dots, c_1, \dots , c_k \} $; \item The sequence of part sizes in each row is weakly decreasing; \item The sequence of subscripts in each row is weakly decreasing; \item Every subscript $1, \dots , k-1$ appears at least once
in the top row; \item If $M_1, M_2, \dots M_{k-2}, M_{k-1}$ are
the largest parts with their respective subscripts in the top row,
then $b_i = d_e \Rightarrow d \in [ M_{e-1}, M_e ]$, setting $M_0
= 1$ and $M_k = c$.
\end{itemize}
\end{definition}

Calling $D_k (n)$ the number of $k$-marked Durfee symbols of $n$,
we then have $D_{k+1} (n) = \eta_{2k} (n)$ (Corollary 13 in
\cite{Andrews}; the $\eta_{2k+1} (n)$ vanish, as do the Atkin-Garvan odd moments). The study of congruence theorems for Durfee
symbols thus informs the study of congruence theorems for standard 
partitions.  To study this relation he defines the
\emph{full rank} of a $k$-marked Durfee symbol:

\begin{definition} Let $\delta$ be a $k$-marked Durfee symbol and let $\tau_i (resp. \, \beta_i)$ be the number of parts in the top (resp. bottom) row with subscript $i$.  Then the $i^{th}$-rank of $\delta$ is
\[ \rho_i ( \delta ) = \left\{
\begin{matrix}
\tau_i - \beta_i - 1    & 1 \leq i < k \\
\tau_i - \beta_i        & i = k
\end{matrix} \right. \, \text{.}
\]
\end{definition}

\begin{definition}The \emph{full rank} of a $k$-marked Durfee symbol $\delta$ is $\rho_1 ( \delta ) + 2 \rho_2 ( \delta ) + 3 \rho_3 ( \delta ) + \dots + k \rho_k ( \delta )$.
\end{definition}

We set $D_k (m_1, \dots , m_k;n)$ to be the number of $k$-marked Durfee symbols with $i^{th}$ rank $m_i$.  In analogy to our previous construction for the rank we call $NF_l (m,n)$ the number of $l$-marked Durfee symbols of $n$ with full rank $m$, and $NF_l (b,c,n)$ the number of $l$-marked Durfee symbols of $n$ with full rank $\equiv b \, (mod \, c)$.

Andrews produces the generating function (Theorems 10 and 7 in
\cite{Andrews}):

\begin{multline}\label{RFunction}
\sum_{n_1, \dots , n_k = -\infty }^{ \infty } \sum_{n \geq 0} D_k (n_1, \dots , n_k;q) {x_1}^{n_1} \dots {x_k}^{n_k} q^n = R_k (x_1 , \dots x_k ; q) \\
= \sum_{i=1}^k \frac{R_1(x_i;q)}{\prod\limits_{j=1 \atop j \neq
i}^k (x_i - x_j) (1 - {x_i}^{-1}{x_j}^{-1})} \, \text{.}
\end{multline}

This theorem in hand, he produces two congruences: that $D_2 (n)
\equiv 0 \, (mod \, 5)$ for $n \equiv 1,4 \, (mod \, 5)$ and $D_3
(n) \equiv 0 \, (mod \, 7)$ for $n \equiv 0,1,5 \, (mod \, 7)$,
because $NF_2 (i,5,n) = NF_2 (j,5,n)$ and $NF_3 (i,7,n) = NF_3
(j,7,n)$ for all $i,j$ in those progressions.  Furthermore, we still have $NF_2
(i,5,n) = NF_2 (j,5,n)$ and $NF_3 (i,7,n) = NF_3 (j,7,n)$ for all
$i,j \neq 0$ in any progression.

As mentioned earlier, the above two theorems are the simplest two
examples of an infinite family of related theorems, which we explore.  We also examine in full detail the behavior of the residue classes for nonprime (odd) modulus.

\section{Proof of Theorem \ref{Main Theorem}}

Let ${\zeta}_c$ be a primitive $c$-th root of unity.  We employ the same basic strategy as Andrews, considering $ \sum_{n=1}^{\infty} \sum_{b=0}^{c-1} NF_l (b,c,n) {{\zeta}_c}^b q^n $.  To prove the general theorem requires the additional observation that, since $N(m,n) = N(-m,n)$, this sum behaves well with respect to sums of conjugate powers of ${\zeta}_c$.  Break the sum down thus:

\begin{multline}\label{ZetaIdent}
\sum_{n=1}^{\infty} \sum_{b=0}^{c-1} NF_l (b,c,n) {{\zeta}_c}^b q^n = R_l ({{\zeta}_c}, {{\zeta}_c}^2, \dots, {{\zeta}_c}^l ; q) \\
 = \sum\limits_{i=1}^l  \frac{ R_1 \left( {{\zeta}_c}^i ; q \right) }{\prod\limits_{j=1 \atop j \neq i }^l \left( {{\zeta}_c}^i - {{\zeta}_c}^j \right) \left(1 - {{\zeta}_c}^{-i-j} \right) }\\ 
 = \sum_{i=1}^l \left( \prod_{j=1 \atop j \neq i }^l \left( {{\zeta}_c}^i - {{\zeta}_c}^j \right) \left(1 - {{\zeta}_c}^{-i-j} \right) \right)^{-1} \cdot \sum_{n=1}^{\infty} \sum_{k=0}^{c-1} \sum_{d=0}^{c-1} {{\zeta}_c}^{i k} N(k,c,cn+d) q^{cn+d} \\
= \sum_{i=1}^l \left( \prod_{j=1 \atop j \neq i }^l \left( {{\zeta}_c}^i - {{\zeta}_c}^j \right) \left(1 - {{\zeta}_c}^{-i-j} \right) \right)^{-1} \cdot \sum_{d=0}^{c-1} q^d \sum_{n \geq 1} \sum_{k=0}^{c-1} {{\zeta}_c}^{i k} N(k,c,cn+d) q^{cn} \quad \text{.}
\end{multline}

Following Atkin, we define $r_{a,b}(q;c;d) = \sum_{n \geq 0} q^n ( N(a,c,cn+d) - N(b,c,cn+d) )$.  Then, for any given $d$, 

\begin{eqnarray*}
\sum_{n \geq 1} N(l,c,cn+d) q^{cn} &= \sum_{n \geq 1} N(l-1,c,cn+d) q^{cn} - r_{l-1,l}(q^c;c;d) \\
                                   &= \sum_{n \geq 1} N(l-2,c,cn+d) q^{cn} - r_{l-2,l}(q^c;c;d) \\
                                   &\vdots \\
                                   &= \sum_{n \geq 1} N(0,c,cn+d) q^{cn} - r_{0,l}(q^c;c;d) \quad \text{.}
\end{eqnarray*}

Now using the evenness of the rank function and the fact that $\sum_{b=0}^{c-1} {{\zeta}_c}^b = 0$, we have

\begin{multline}
\left( \sum_{n \geq 1} N(0,c,cn+d) q^{cn} - r_{0,l} (q^c;c;d) \right) {{\zeta}_c}^{i \cdot 0} \\
+ \left( \sum_{n \geq 1} N(1,c,cn+d) q^{cn} - r_{1,l} (q^c;c;d) \right) {{\zeta}_c}^{i \cdot 1} + \dots \\
+ \left( \sum_{n \geq 1} N(l,c,cn+d) q^{cn} \right) {{\zeta}_c}^{i \cdot l} + \left( \sum_{n \geq 1} N(l+1,c,cn+d) q^{cn} \right) {{\zeta}_c}^{i \cdot (l+1)} \\
+ \dots + \left( \sum_{n \geq 1} N(c-1,c,cn+d) q^{cn} - r_{1,l} (q^c;c;d) \right) {{\zeta}_c}^{i \cdot (c-1)} = 0 \, \text{.}
\end{multline}

\noindent (For use later we note that it matters in the above calculation that $i \not\equiv 0 \ \emph{mod} \ c$ in this context, but its value otherwise is irrelevant; if $c$ is nonprime and $gcd(i,c) \neq 1$, we have merely employed the same identity $\frac{c}{gcd(i,c)}$ times.)

Thus, gathering the $N(k,c,cn+d)$ terms and recalling that $c=2l+1$, 
\[
\sum_{k=0}^{c-1} {{\zeta}_c}^{ik} \sum_{n \geq 1} N(k,c,cn+d) q^{cn} = r_{0,l} (q^c ;c;d) + \sum_{g=1}^{l-1} r_{g,l} (q^c ;c;d) \left( {{\zeta}_c}^{ig} + {{\zeta}_c}^{i(-g)} \right) \ \text{.}
\]

Hence
\begin{multline}
\sum_{n=1}^{\infty} \sum_{b=0}^{c-1} NF_l (b,c,n) {{\zeta}_c}^b q^n = \sum_{i=1}^l \left( \prod_{j=1 \atop j \neq i }^l \left( {{\zeta}_c}^i - {{\zeta}_c}^j \right) \left(1 - {{\zeta}_c}^{-i-j} \right) \right)^{-1} \\
\cdot \sum_{d=0}^{c-1} q^d \left( r_{0,l} (q^c ;c;d) + \sum_{g=1}^{l-1} r_{g,l} (q^c ;c;d) \left( {{\zeta}_c}^{ig} + {{\zeta}_c}^{i(-g)} \right) \right) \ \text{.}
\end{multline}

For any $n$, then, we have by equation of coefficients in powers of $q$ that \vspace{-0.1in} 
\begin{multline}\label{EqnCoeff}
\sum_{b=0}^{c-1} NF_l (b,c,n) {{\zeta}_c}^b = \sum_{i=1}^l \left( \prod_{j=1 \atop j \neq i }^l \left( {{\zeta}_c}^i - {{\zeta}_c}^j \right) \left(1 - {{\zeta}_c}^{-i-j} \right) \right)^{-1} \\
\cdot \left( N(0,c,n) - N(l,c,n) + \sum_{g=1}^{l-1} \left( N(g,c,n) - N(l,c,n) \right) \left( {{\zeta}_c}^{ig} + {{\zeta}_c}^{i(-g)} \right) \right) \  \text{.}
\end{multline}

To prove the theorem, we must first show that the right-hand side of (\ref{EqnCoeff}) is an integer.  The constant term that appears before the sum contributes 0: notice that $\left( \prod\limits_{j=1 \atop j \neq i }^l \left( {{\zeta}_c}^i - {{\zeta}_c}^j \right) \left( 1 - {{\zeta}_c}^{-i-j} \right) \right) = {{\zeta}_c}^{i(l-1)} \left( \prod\limits_{j=1 \atop j \neq i }^l \left( 1 - {{\zeta}_c}^{-i+j} \right) \left(1 - {{\zeta}_c}^{-i-j} \right) \right) $, and the exponents $ \{ -i+j, -i-j \, \vert \, 1 \leq j \leq l, j \neq i \} $ are precisely $ \{ 1, \dots , c-1 \} \setminus \{ 0, i, 2i \} $ when reduced \emph{mod} $c$.  Since $\prod_{i=1}^{c-1} \left( 1 - {{\zeta}_c}^i \right) = c$, we simplify thus: 

\begin{multline}
\left( N(0,c,n) - N(l,c,n) \right) \cdot \sum_{i=1}^l \left( \prod_{j=1 \atop j \neq i }^l \left( {{\zeta}_c}^i - {{\zeta}_c}^j \right) \left(1 - {{\zeta}_c}^{-i-j} \right) \right)^{-1} \\ 
\shoveleft{ \quad = \left( N(0,c,n) - N(l,c,n) \right) \cdot \frac{1}{c} \cdot \sum_{i=1}^l {{\zeta}_c}^{-i(l-1)} \left( 1 - {{\zeta}_c}^{-2i} \right) \left( 1 - {{\zeta}_c}^{-i} \right) } \\ 
\shoveleft{ \quad = \left( N(0,c,n) - N(l,c,n) \right) \cdot \frac{1}{c} \cdot \sum_{i=1}^l \left( {{\zeta}_c}^{-i(l-1)} + {{\zeta}_c}^{-i(l+2)} - {{\zeta}_c}^{-i(l+1)} - {{\zeta}_c}^{-i(l)} \right) } \\
\shoveleft{ \quad = \left( N(0,c,n) - N(l,c,n) \right) \cdot \frac{1}{c} \cdot \sum_{i=1}^l \left( {{\zeta}_c}^{-i(l-1)} + {{\zeta}_c}^{i(l-1)} - {{\zeta}_c}^{-i(l+1)} - {{\zeta}_c}^{i(l+1)} \right) } \\
\shoveleft{ \quad = \left( N(0,c,n) - N(l,c,n) \right) \cdot \frac{1}{c} \cdot \sum_{i=1}^{2l} \left( {{\zeta}_c}^{i(l-1)} - {{\zeta}_c}^{-i(l+1)} \right) } \\
\shoveleft{ \quad = \left( N(0,c,n) - N(l,c,n) \right) \cdot \frac{1}{c} \cdot (-1 -(-1)) = 0 \ \text{.} } \hfil
\end{multline}

There remains the second term, which contributes a nonzero integer:
\begin{multline}
\sum_{i=1}^l \left( \prod_{j=1 \atop j \neq i }^l \left( {{\zeta}_c}^i - {{\zeta}_c}^j \right) \left(1 - {{\zeta}_c}^{-i-j} \right) \right)^{-1} \cdot \sum_{g=1}^{l-1} \left( N(g,c,n) - N(l,c,n) \right) \left( {{\zeta}_c}^{ig} + {{\zeta}_c}^{i(-g)} \right) \\
= \sum_{i=1}^l {{\zeta}_c}^{-i(l-1)} \left( 1 - {{\zeta}_c}^{-2i} \right) \left( 1 - {{\zeta}_c}^{-i} \right) \cdot \frac{1}{c} \cdot \sum_{g=1}^{l-1} \left( N(g,c,n) - N(l,c,n) \right) \left( {{\zeta}_c}^{ig} + {{\zeta}_c}^{i(-g)} \right) \\
 = \frac{1}{c} \cdot \sum_{i=1}^l \sum_{g=1}^{l-1} \left( N(g,c,n) - N(l,c,n) \right) \left[ {{\zeta}_c}^{-i(l-g-1)} + {{\zeta}_c}^{-i(l+g-1)} + {{\zeta}_c}^{-i(l-g+2)} \right. \\ 
 \left. \qquad + {{\zeta}_c}^{-i(l+g+2)} - {{\zeta}_c}^{-i(l-g)} - {{\zeta}_c}^{-i(l+g)} - {{\zeta}_c}^{-i(l-g+1)} - {{\zeta}_c}^{-i(l+g+1)} \right] \\
 = \frac{1}{c} \cdot \sum_{g=1}^{l-1} \left( N(g,c,n) - N(l,c,n) \right) \sum_{i=1}^l \left[ {{\zeta}_c}^{-i(l-g-1)} + {{\zeta}_c}^{i(l-g+2)} + {{\zeta}_c}^{-i(l-g+2)} \right. \\
 \left. \qquad + {{\zeta}_c}^{i(l-g-1)} - {{\zeta}_c}^{-i(l-g)} - {{\zeta}_c}^{i(l-g+1)} - {{\zeta}_c}^{-i(l-g+1)} - {{\zeta}_c}^{i(l-g)} \right] \\
 = \frac{1}{c} \cdot \sum_{g=1}^{l-1} \left( N(g,c,n) - N(l,c,n) \right) \sum_{i=1}^{c-1} \left[ {{\zeta}_c}^{i(l-g+2)} + {{\zeta}_c}^{i(l-g-1)} - {{\zeta}_c}^{i(l-g+1)} - {{\zeta}_c}^{i(l-g)} \right] \\
 = \frac{1}{c} \cdot \sum_{g=1}^{l-1} \left( N(g,c,n) - N(l,c,n) \right) \cdot \epsilon \quad \text{,}
\end{multline}

\noindent where $\epsilon = 0$ if $g \neq l-1$ and $\epsilon = c$ if $g = l-1$.

Thus, the right-hand side of (\ref{EqnCoeff}) is an integer, and so (\ref{EqnCoeff}) is a polynomial of degree $c-1$ in ${{\zeta}_c}$ over the integers.  We can particularly evaluate

\begin{equation}\label{IntegerForm}
\sum_{b=0}^{c-1} NF_l (b,c,n) {{\zeta}_c}^b = N(l-1,c,n) - N(l,c,n) \quad \text{.}
\end{equation}

\noindent When $c$ is prime, the irreducibility of the minimal polynomial $1+x+ \dots + x^{c-1}$ for $\zeta_c$ leads immediately to equality of the coefficients for $b \neq 0$ (if the coefficients are unequal, subtract from (\ref{IntegerForm}) the equation $NF_l (c-1,c,n) (1+ \zeta_c + \dots + {\zeta_c}^{c-1}) = 0$ to obtain a new integer polynomial in $\zeta_c$ of lower degree, contradicting minimality).  

To prove the theorem in the nonprime case, we must recall the symmetries of the original setting, identity (\ref{RFunction}).  We evaluated this equation at $x_i = {{\zeta}_c}^i$ to obtain (\ref{ZetaIdent}), which is an identity of elements in $\mathbb{Q}({\zeta}_c)[[q]]$, power series with coefficients in $\mathbb{Q}({\zeta}_c)$.  But elements of this field can be represented non-uniquely by polynomials in ${\zeta}_c$ of degree less than $c$ when $c$ is nonprime: $0 = {{\zeta}_9}^1 + {{\zeta}_9}^4 + {{\zeta}_9}^7$ and $0 = {{\zeta}_9}^2 + {{\zeta}_9}^5 + {{\zeta}_9}^8$.  So we take an intermediate step: evaluate (\ref{RFunction}) at $x_i = z^i$ to obtain

\begin{equation}\label{ZForm}
\sum_{n=1}^{\infty} \sum_{b > -\infty}^{\infty} NF_l (b,n) z^b q^n = \sum\limits_{i=1}^l  \frac{ R_1 \left( z^i ; q \right) }{\prod\limits_{j=1 \atop j \neq i }^l \left( z^i - z^j \right) \left(1 - z^{-i-j} \right) }
\end{equation}

\noindent which is an identity of elements in $\mathbb{Q}[[q]]((z))$.  The coefficient of any particular $q^n$ is a finite symmetric Laurent polynomial in $z$.  Taking the quotient in this ring by the ideal $(z^p - 1)$, we obtain the identity

\begin{equation}\label{ZModZP}
\sum_{n=1}^{\infty} \sum_{b =0}^{c-1} NF_l (b,c,n) z^b q^n \equiv \sum\limits_{i=1}^l  \frac{ R_1 \left( z^i ; q \right) }{\prod\limits_{j=1 \atop j \neq i }^l \left( z^i - z^j \right) \left(1 - z^{-i-j} \right) }
\end{equation}

\noindent where the equivalence is one of cosets in $\mathbb{Q}[[q]]((z)) / (z^p - 1)$.  There is a unique representative of any such coset in which the coefficients on $q^n$ are true polynomials in $z$ of degree less than $c$.  These polynomials have exactly the coefficients given on the left hand side of the equivalence.  To complete the proof, we note that the map $z \rightarrow z^m$ simply permutes terms of the sum on the right-hand side, if and only if $gcd(m,c)=1$.  Some such $m$ will map any $z^a \rightarrow z^b$ for any given pair of $a$ and $b$ with $gcd(a,c) = gcd(b,c)$.  But permuting the terms of a sum with a finite number of nonzero terms does not alter the result, so the right-hand side is fixed under this mapping.  Because the representatives of degree less than $c$ are unique, the left-hand side must be fixed as well.  In particular, $NF_l (a,c,n) = NF_l (b,c,n)$ as long as $gcd(a,c) = gcd(b,c)$. $\Box$

If we know something about the difference $N(l-1,c,cn+d) - N(l,c,cn+d)$, we can now say something about the behavior of the $l$-ranks.  Work of Atkin and Swinnerton-Dyer \cite{ASD} yields the arithmetic progressions mentioned by Andrews, for $c=5$ and $c=7$, in which the difference is identically 0 and equidistribution of the $l$-ranks is achieved.  Study of the difference for additional prime $c$ has been made by Atkin and collaborators Hussain \cite{AH} and O'Brien \cite{OBrien}: specifically $c=11, 13, 17, \, \text{and} \, 19$.  Results on the moments of ranks which can inform use of these theorems, especially from the viewpoint of automorphic forms, can also be found in the work of Stephanie Treneer with Scott Ahlgren \cite{AhlgrenTreneer}, Bringmann as previously cited (\cite{BGM} and \cite{BK}), and additionally in work with Ken Ono and R. C. Rhoades \cite{BrOnRh}.  In addition, congruences for the standard rank moments were studied by Garvan \cite{Garvan} using a connection to the moments of the $\emph{crank}$ statistic and relations on the coefficients of half-integer weight Hecke eigenforms; those techniques could also be useful in finding explicit congruences for Durfee symbols.

\section{Nonprime Moduli}

We now turn to a deeper examination of nonprime $c$.  No longer is the polynomial $1+x+x^2 + \dots + x^{c-1}$ irreducible over the integers, so the populations of the various divisor-groups of residue classes mod $c$ are no longer necessarily equal.  However, if we can find $N(0,c,n) - N(d,c,n)$ for all $d \, \vert \, c$, we can state a congruence theorem for $D_l (n)$ modulo $c$.

We do this by observing the behavior of $R_l ({{\zeta}_c}^d, {{\zeta}_c}^{2d}, \dots , {{\zeta}_c}^{ld} )$.  Assigning $x_i = {{\zeta}_c}^{di}$ in Equation (\ref{RFunction}) (Theorem 7 of \cite{Andrews}) and simplifying using standard properties of cyclotomic polynomials, we have 

\begin{equation} R_l ({{\zeta}_c}^d, {{\zeta}_c}^{2d}, \dots , {{\zeta}_c}^{ld} ) = \sum_{n = 0}^{\infty} q^n  \sum_{r | c} NF_l(r,c,n) \mu \left( \frac{c}{gd} \right) \frac{\phi \left( \frac{c}{r} \right) }{\phi \left( \frac{c}{gd} \right)}
\end{equation}

\noindent with $\mu$ the standard M\"{o}bius function, $\phi$ the totient, and $g=gcd(r,c/d)$.

The coefficients involved are relatively small.  By way of example we use later, \[ R_4 ({{\zeta}_9}, {{\zeta}_9}^2 , {{\zeta}_9}^3 , {{\zeta}_9}^4 ; q) = \sum_{n \geq 0} q^n \left( NF_4 (0,9,n) + 0 \cdot NF_4(1,9,n) - NF_4 (3,9,n) \right) \, \text{,} \] \noindent and \[ R_4 ({{\zeta}_9}^3, {{\zeta}_9}^6 , {{\zeta}_9}^9 , {{\zeta}_9}^{12} ; q) = \sum_{n \geq 0} q^n \left( NF_4 (0,9,n) - 3 NF_4 (1,9,n) + 2 NF_4 (3,9,n) \right) \, \text{.} \]

Calculating this value for each $d$ strictly dividing $c$ gives us a system of $d(c)-1$ linear equations in the $NF_l(d,c,n)$ (where $d(c)$ is the divisor function) that we can solve explicitly for the differences $NF_l(0,c,n) - NF_l(d,c,n)$.

At first glance, assigning $x_i = {{\zeta}_c}^{di}$ in Equation (\ref{RFunction}), where no longer $d=1$ as in the main theorem, produces singularities in the terms $\frac{1}{(x_i - x_j)(1 - x_i^{-1}x_j^{-1})}$ when $j \equiv \pm i \, (mod \, c/d )$.  These singularities are, of course, removable by repeated application of L'H\^{o}pital's rule.

The case $c=9$ is the first opportunity to employ the method, the most tractable to calculate explicitly for illustrative purposes, and an interesting example in its own right.  Begin with the $l=4$ case of (\ref{RFunction}):

\[
R_4 (x_1, x_2, x_3, x_4;q) = \sum_{i=1}^4 \frac{R_1 (x_i; q)}{\prod\limits_{j=1 \atop j \neq i}^4 (x_i - x_j)(1-x_i^{-1}x_j^{-1})} \, \text{.}
\]

We know that

\begin{multline}\label{R4Basic}
R_4 ({{\zeta}_9}, {{\zeta}_9}^2 , {{\zeta}_9}^3 , {{\zeta}_9}^4 ; q) = \sum_{n \geq 0} q^n \left( NF_4 (0,9,n) - NF_4 (3,9,n) \right) \\
\quad = \sum_{n \geq 0} q^n \left( N(3,9,n) - N(4,9,n) \right) \, \text{.}
\end{multline}

Already we can state an interesting congruence: a conjecture of Richard Lewis \cite{Lewis} proved by Nicholas Santa Gadea \cite{SG} states that $N(3,9,3n) = N(4,9,3n)$.  Thus $NF_4 (0,9,3n) = NF_4 (3,9,3n) = NF_4 (6,9,3n)$ and, since $NF_4(i,9,n) = NF_4(j,9,n)$ for the 6 residue classes $3 \nmid i,j$, we have proved Theorem \ref{FourMarkedModThree}.

To say more regarding the behavior of $D_4$ mod 9, we need to know the difference $NF_4 (0,9,n) - NF_4 (1,9,n)$.  To obtain this we calculate, for $d=3$, 

\begin{multline}
R_4 ({{\zeta}_9}^3, {{\zeta}_9}^6 , {{\zeta}_9}^9 , {{\zeta}_9}^3 ; q) = R_4 ({\zeta}_3, {{\zeta}_3}^2 , 1 , {\zeta}_3 ; q) \\
\quad = \sum_{n \geq 0} q^n \left( NF_4 (0,9,n) - 3 NF_4 (1,9,n) + 2 NF_4 (3,9,n) \right)
\end{multline}

\noindent in terms of $R_1 (\zeta_3 ; q)$.

The strategy is to replace, one by one, each of the $x_i$ by functions of $x_1$ which replicate the relations of the ${{\zeta}_3}^i$: $x_4$ by $x_1$, $x_3$ by 1, and $x_2$ by $x_1^{-1}$.  At each step we obtain a small number of singularities we can remove.  First, replace $x_4$ by $x_1$.

\begin{multline}
R_4 (x_1, x_2, x_3, x_1;q) = {\lim_{x_4 \rightarrow x_1}} R_4 (x_1, x_2, x_3, x_4;q) \\
=  {\lim_{x_4 \rightarrow x_1}} \left( \frac{R_1 (x_1;q)}{\prod_{2 \leq j \leq 4} (x_1 - x_j)(1 - x_1^{-1} x_j^{-1})} + \frac{R_1 (x_4;q)}{\prod_{1 \leq j \leq 3}(x_4 - x_j)(1 - x_4^{-1}x_j^{-1})} \right) \\
+ \frac{R_1 (x_2;q)}{{(x_2 - x_1)}^2 (x_2 - x_3){(1 - x_2^{-1} x_1^{-1})}^2 (1 - x_2^{-1} x_3^{-1})} \\
+ \frac{R_1 (x_3;q)}{{(x_3 - x_1)}^2 (x_3 - x_2){(1 - x_3^{-1}x_1^{-1})}^2 (1 - x_3^{-1}x_2^{-1})} \\
=  \frac{R_1 (x_2;q)}{{(x_2 - x_1)}^2 (x_2 - x_3){(1 - x_2^{-1} x_1^{-1})}^2 (1 - x_2^{-1} x_3^{-1})} \hfill \\
+ \frac{R_1 (x_3;q)}{{(x_3 - x_1)}^2 (x_3 - x_2){(1 - x_3^{-1}x_1^{-1})}^2 (1 - x_3^{-1}x_2^{-1})} \\
+ {\lim_{x_4 \rightarrow x_1}} \left( \frac{1}{(x_4 - x_1)(1-x_4^{-1}x_1^{-1})} \frac{1}{\prod\limits_{i=1,4 \atop j=2,3} (x_i - x_j) (1 - x_i^{-1} x_j^{-1})} \right. \\
\times \left( R_1 (x_4;q) (x_1 - x_2)(x_1 - x_3)(1 - x_1^{-1} x_2^{-1})(1 - x_1^{-1} x_3^{-1}) - \right. \\
 \left. \left. R_1 (x_1;q) (x_4 - x_2)(x_4 - x_3)(1 - x_4^{-1}x_2^{-1})(1 - x_4^{-1}x_3^{-1}) \right) \right) \text{.}
\end{multline}

After differentiation and taking the limit, we obtain

\begin{multline}
R_4 (x_1, x_2, x_3, x_1;q) = \frac{R_1 (x_2;q)}{{(x_2 - x_1)}^2 (x_2 - x_3){(1 - x_2^{-1} x_1^{-1})}^2 (1 - x_2^{-1} x_3^{-1})} \\
+ \frac{R_1 (x_3;q)}{{(x_3 - x_1)}^2 (x_3 - x_2){(1 - x_3^{-1}x_1^{-1})}^2 (1 - x_3^{-1}x_2^{-1})} \\
+ \frac{\frac{\partial}{\partial x_1} R_1 (x_1; q)}{(x_1 - x_2) (x_1 - x_3) (1 - x_1^{-1}x_2^{-1})(1 - x_1^{-1} x_3^{-1}) (1 - x_1^{-2})} \\
- \frac{R_1 (x_1; q)(\frac{1}{x_1 - x_2} + \frac{1}{x_1 - x_3} + \frac{x_1^{-2}x_2^{-1}}{1-x_1^{-1}x_2^{-1}} + \frac{x_1^{-2}x_3^{-1}}{1-x_1^{-1}x_3^{-1}})}{(x_1 - x_2) (x_1 - x_3) (1 - x_1^{-1}x_2^{-1})(1 - x_1^{-1} x_3^{-1}) (1 - x_1^{-2})} \, \text{.}
\end{multline}

For the next step we replace $x_3$ by 1.  In the case of $d=3$, replacing $x_{c/3}$ by 1 produces no singularities, and so we need not differentiate.  (For any $c$, 3 is the only divisor where this degeneracy ever occurs; for any other potential divisor of $c$, $\lfloor \frac{l}{c/d} \rfloor >2$ since $c=2l+1$, so this replacement step would produce singularities in the denominator factors $(x_{\frac{kc}{d}} - x_{\frac{hc}{d}})$ and $(1 - x_{\frac{kc}{d}}^{-1}x_{\frac{hc}{d}}^{-1})$ in Equation (\ref{RFunction}).)  For $c=9$, we obtain

\begin{multline}
R_4 (x_1, x_2, 1, x_1;q) = \frac{R_1 (x_2;q)}{{(x_2 - x_1)}^2 (x_2 - 1){(1 - x_2^{-1} x_1^{-1})}^2 (1 - x_2^{-1})} \\
+ \frac{R_1 (1;q)}{{(1 - x_1)}^2 (1 - x_2){(1 - x_1^{-1})}^2 (1 - x_2^{-1})} \\
+ \frac{\frac{\partial}{\partial x_1} R_1 (x_1; q)}{(x_1 - x_2) (x_1 - 1) (1 - x_1^{-1}x_2^{-1})(1 - x_1^{-1}) (1 - x_1^{-2})} \\
- \frac{R_1 (x_1; q)(\frac{1}{x_1 - x_2} + \frac{1}{x_1 - 1} + \frac{x_1^{-2}x_2^{-1}}{1-x_1^{-1}x_2^{-1}} + \frac{x_1^{-2}}{1-x_1^{-1}})}{(x_1 - x_2) (x_1 - 1) (1 - x_1^{-1}x_2^{-1})(1 - x_1^{-1}) (1 - x_1^{-2})} \, \text{.}
\end{multline}

It remains to replace $x_2$ by $x_1^{-1}$.

\begin{multline}
R_4 (x_1, x_1^{-1}, 1, x_1; q) = {\lim_{x_2 \rightarrow x_1^{-1}}} R(x_1, x_2, 1, x_1;q) = \frac{R_1(1;q)}{{(1-x_1)}^3{(1-x_1^{-1})}^3} \\
+ {\lim_{x_2 \rightarrow x_1^{-1}}} \left( \frac{1}{{(x_1-x_2)}^2 (x_1 - 1) {(1 - x_1^{-1}x_2^{-1})}^2 (1-x_1^{-1})(1-x_2^{-1})(1-x_1^{-2})(x_2 - 1)} \right. \\
\times \left( R_1 (x_2;q)((x_1-1)(1-x_1^{-1})(1-x_1^{-2})) + ((x_1 - x_2)(1-x_1^{-1}x_2^{-1})(x_2 - 1)(1-x_2^{-1}) ) \right .\\
\times \left. \left. \left( \frac{\partial}{\partial x_1} R_1 (x_1; q) - R_1(x_1; q) (\frac{1}{x_1 - x_2} + \frac{1}{x_1 - 1} + \frac{x_1^{-2}x_2^{-1}}{1-x_1^{-1}x_2^{-1}} + \frac{x_1^{-2}}{1-x_1^{-1}}) \right) \right) \right) \, \text{.}
\end{multline}

We differentiate (twice) with respect to $x_2$.  Using the identity 

\[
{\lim_{x_2 \rightarrow x_1^{-1}} } \frac{\partial^2}{\partial {x_2}^2} R_1 (x_2;q) = {x_1}^4 \frac{\partial^2}{\partial {x_1}^2} R_1 (x_1;q) + 2 {x_1}^3 \frac{\partial}{\partial x_1} R_1 (x_1;q) \, \text{,}
\]

\noindent which follows from the symmetry of $R_1 (x;q)$ that $R_1(x, q) = R_1(x^{-1}, q)$, we obtain in the limit

\begin{multline}
R_4 (x_1, x_1^{-1}, 1, x_1; q) = \frac{R_1(1;q)}{{(1-x_1)}^3{(1-x_1^{-1})}^3} + \frac{-{x_1}^{-4}}{2{(1-x_1)}^3{(1-{x_1}^{-1})}^3{(1-{x_1}^{-2})}^3} \\
 \times \left[ \frac{\partial^2}{\partial {x_1}^2} R_1 (x_1;q) {x_1}^4 {(x_1 - 1)}^2 {(1-{x_1}^{-1})}^2 (1-{x_1}^{-2}) \right. \\
+ 2 \frac{\partial}{\partial x_1} R_1 (x_1;q) {(1-x_1 )}^2 {(1-x_1^{-1})}^2 (-{x_1}^3 - 2 {x_1}^2 - 2 x_1 ) \\
+ \left. 2 R_1 (x_1 ; q) {(1-{x_1}^2)}^2 (1-{x_1}^{-2}) \right] \, \text{.} 
\end{multline}

We have now removed all the troublesome singularities and can set in the last identity $x_1 = {\zeta}_3$ to evaluate 

\begin{multline}
R_4 ({\zeta}_3, {{\zeta}_3}^{2}, 1, {\zeta}_3; q) = \sum_{n \geq 0} q^n ( NF_4 (0,9,n) - 3 NF_4 (1,9,n) + 2 NF_4 (3,9,n) ) \\
= \frac{R_1(1;q)}{27} + \frac{-{{\zeta}_3}^2}{54 {(1-{\zeta}_3)}^3} \hfil \\
\quad \times \left[ 9 ({\zeta}_3 - {{\zeta}_3}^2 ) \frac{\partial^2}{\partial {z}^2} R_1 (z;q) {\Big|}_{z={\zeta}_3} + 18 \frac{\partial}{\partial z} R_1 (z;q) {\Big|}_{z={\zeta}_3} + 6 (1 - {{\zeta}_3}^2 ) R_1 ({\zeta}_3;q) \right] \\
 = \frac{1}{54} \left( 2 R_1 (1;q) + 3 {{\zeta}_3}^2 \frac{\partial^2}{\partial {z}^2} R_1 (z;q) {\Big|}_{z={\zeta}_3} + 2 ({\zeta}_3 - 1 ) \frac{\partial}{\partial z} R_1 (z;q) {\Big|}_{z={\zeta}_3} - 2 R_1 ({\zeta}_3;q) \right) \text{.}
\end{multline}

We wish to rewrite this formula in terms of the rank classes $N(j,n)$.  The termwise first and second derivatives of $N(j,n) z^j q^n$, $j N(j,n) z^{j-1} q^n$ and $j(j-1) N(j,n) z^{j-2} q^n$ respectively, group themselves thus by the residue class of $j$ modulo 3 when evaluated at $z=\zeta_3$:

\begin{multline}\label{ZetaThreeSum}
R_4 ({\zeta}_3, {{\zeta}_3}^{2}, 1, {\zeta}_3; q) \\
= \frac{1}{54} \sum_{n \geq 0} q^n \sum_{k=-\infty}^{\infty} \left( (27 k^2 -3k) N(3k,n) - 6k N(3k+1,n) + 2 N(3k+2,n) \right. \\
\quad +  \zeta_3 \left( (27k^2 + 15k) N(3k+1,n) - (6k+4) N(3k+2,n) \right) \\
\qquad +\left. {\zeta_3}^2 \left( (27k^2+33k+8)N(3k+2,n) -6k N(3k,n) \right) \right) \\
= \frac{1}{54} \sum_{n \geq 0} q^n \sum_{k=-\infty}^{\infty} (27k^2 + 3k) N(3k,n) - 6k N(3k+1,n) - (27k^2+33k+6) N(3k+2,n)  \\ 
+ \zeta_3 \left( 6k N(3k,n) + (27k^2 + 15k) N(3k+1,n) - (27k^2 + 39k + 12) N(3k+2,n) \right) \, \text{.}
\end{multline}

The sums are finite, since the rank $N(k,n)$ is identically 0 for $|k| > n-1$.  We can simplify the sum above by recalling that, due to the evenness of the rank function, for any $j$

\[ \sum_{k=-j}^j k N(3k,n) = \sum_{k=-j}^{j-1} k N(3k+1,n) + (k+1) N(3k+2,n) = 0  \]

and 

\[ \sum_{k=-j}^j k^2 N(3k+1,n) = \sum_{k=-j}^{j-1} (k+1)^2 N(3k+2,n) \, \text{.} \]

With these two identities the $\zeta_3$ term of (\ref{ZetaThreeSum}) wholly vanishes.  (We knew it must, since of course $R_4 (\zeta_3 , {\zeta_3}^2, 1, \zeta_3 ; q)$ has integral coefficients.)

Upon discarding the vanishing $\zeta_3$ term and simplifying the remainder with the relations above we have that 

\begin{multline}
NF_4 (0,9,n) - 3 NF_4 (1, 9, n) + 2 NF_4 (3,9,n) \\
= \sum_{k=-\infty}^{\infty} \left( \frac{k^2}{2} N(3k,n) - \frac{k(k+1)}{2} N(3k+2,n) \right) \\
= \sum_{k=1}^{\infty} k^2 N(3k,n) - \frac{k(k+1)}{2} \left( N(3k+1,n) + N(3k+2,n) \right) \, \text{.}
\end{multline}

Thus, since we already know from Equation (\ref{R4Basic}) that $NF_4 (0, 9, n) - NF_4 (3,9,n) = N(3,9,n) - N(4,9,n)$,

\begin{multline}
NF_4 (0,9,n) - NF_4 (1,9,n) = - \frac{2}{3} \left( N(3,9,n) - N(4,9,n) \right) \\
+ \frac{1}{3} \sum_{k=1}^{\infty} k^2 N(3k,n) - \frac{k(k+1)}{2} \left( N(3k+1,n) + N(3k+2,n) \right) \, \text{.}
\end{multline}

Putting these all together, we have

\begin{multline} D_4 (n) = \sum_{i=0}^{8} NF_4 (i,9,n) = NF_4 (0,9,n) + 6 NF_4 (1,9,n) + 2 NF_4 (3,9,n) \\ = NF_4 (0,9,n) + 6 (NF_4 (0,9,n) - (NF_4 (0,9,n) - NF_4 (1,9,n))) \\
+ 2 (NF_4 (0,9,n) - (NF_4 (0,9,n) - NF_4 (3,9,n))) \\
 = 9 NF_4 (0,9,n) + 2 (N(3,9,n) - N(4,9,n)) \hfil \\
 - \sum_{k=1}^{\infty} 2k^2 N(3k,n) - k(k+1) (N(3k+1,n) + N(3k+2,n)) \\
\equiv 2 (N(3,9,n) - N(4,9,n)) \hfil \\
 - \sum_{k=1}^{\infty} 2k^2 N(3k,n) - k(k+1) (N(3k+1,n) + N(3k+2,n)) \, (mod \, 9) \, \text{.}
\end{multline}

For $n \equiv 0, 1, 2 \, (mod \, 3)$, the identities of \cite{SG} provide specializations of this identity when we dissect the sum over $k$ by the residue classes of $k$ modulo 3.

In the case of general $c$ and divisor $d$, we perform variable replacements patterned on those we saw above.  We replace $x_{i+\frac{kc}{d}}$ with $x_i$ for $0 < i < \frac{c}{d}$, replace $x_{\frac{kc}{d}}$ with 1, and finally replace $x_{\frac{c}{d} - i}$ with $x_i^{-1}$ for $0 < i \leq \lfloor \frac{c}{2d} \rfloor $.  We eventually encounter derivatives of order up to $2d$, in order to clear singularities.  When we then evaluate Theorem 7 at $x_i = {{\zeta}_c}^{di}$, a great deal of simplification can occur by working with the evenness of the rank function.  The process is straightforward and might even be automatable.

\section{Acknowledgments}

The content of this article was primarily produced as part of the author's thesis work, under the guidance of George Andrews as the author's advisor.  His support and mentorship were invaluable at every stage of production.

Careful refereeing work improved much exposition and occasionally argument in this paper, for which the author is grateful.  The referee was clearly well-informed as to the current state of the field and provided several of the pointers to current work mentioned above, which the reader will hopefully find useful.

\end{document}